\documentclass[12pt]{elsarticle}

\usepackage{fullpage}

\begin{document}
 
\begin{frontmatter}

 

\title{Delay Terms in the Slow Flow}

\author[a]{Si Mohamed Sah\corref{cor1}}
 \ead{smsah@kth.se}

 \author[b]{ Richard H. Rand}
 
 \cortext[cor1]{Corresponding author}

  \address[a]{Nanostructure Physics, KTH Royal Institute of Technology, \\
Roslagstullsbacken 21, SE-106 91, Stockholm, Sweden.}
 
 \address[b]{Dept. of Mathematics, Dept. of Mechanical \& Aerospace Engineering, Cornell University.\\
 Ithaca, NY 14853, USA.}

\begin{abstract}

This work concerns the dynamics of nonlinear systems that are subjected to delayed self-feedback. Perturbation methods applied to such systems give rise to slow flows which characteristically contain delayed variables. 
We consider two approaches to analyzing Hopf bifurcations in such slow flows.  
In one approach, which we refer to as approach \bf I\rm,  we follow many researchers in replacing the delayed variables in the slow flow with non-delayed variables, thereby reducing the DDE slow flow to an ODE. In a second approach, which we refer to as approach \bf II\rm, we keep the delayed variables in the slow flow.  By comparing these two approaches we are able to assess the accuracy of making the simplifying assumption which replaces the DDE slow flow by an ODE.  We apply this comparison to two examples, Duffing and van der Pol equations with self-feedback.
 
\end{abstract}

\begin{keyword}

Slow flow, Delay, Duffing, Van der Pol, Hopf bifurcation.
\end{keyword}

\end{frontmatter}

\section{INTRODUCTION}

It is known that ordinary differential equations (ODEs) are used as models to better understand phenomenon occurring in biology, physics and engineering. Although these models present a good approximation of the observed phenomenon, in many cases they fail to capture the rich dynamics observed in natural or technological systems. Another approach which has gained interest in modeling systems is  the inclusion of time delay terms in the differential equations resulting in delay-differential equations (DDEs).  DDE's have found application in many systems, including rotating machine tool vibrations \cite{Moon}, gene copying dynamics \cite{Verdugo},  laser dynamics \cite{Wirkus} and many other examples.\\

Despite their simple appearance,  delay-differential equations (DDEs) have several features that make their analysis a challenging task. For example, when investigating a delay-differential equation (DDE) by use of a perturbation method,
one is often confronted with a slow flow which contains delay terms.  It is usually argued
that since the parameter of perturbation, call it $\epsilon$, is small, $\epsilon <<1$, the delay terms 
which appear in the slow flow may be replaced by the same term without delay, see e.g.\,\cite{Ji, Maccari, Hu, Wahi, Atay, Suchorsky, Wirkus, Morrison}.
The purpose of the present paper is to analyze the slow flow with the delay terms left in it, and to compare the resulting
approximation with the usual one in which the delay terms have been replaced by terms without delay.\\

The general class of DDEs that we are interested in is of the form
\begin{equation}
\ddot{x}+x=\epsilon f(x,x_d)
\end{equation}
where $x_d=x(t-T)$, where $T= $ delay.\\

As an example we choose the Duffing equation with delayed self-feedback.
\begin{equation}
\ddot{x}+x=\epsilon ~\left[-\alpha \dot x -\gamma x^3+ k ~x_d\right]
\label{duf}
\end{equation}

The situation here is that when there is no feedback ($k=0$), the 
Duffing equation does not exhibit a limit cycle.
However it turns out that for $k>\alpha$ a stable limit cycle is born in a Hopf bifurcation 
for a critical value of delay $T$ that depends on $k$.
Further increases in $T$ produce another Hopf, which sees the stable limit cycle disappear.

See Figure 1 which shows a plot of the Hopfs in the $k-T$ parameter plane, obtained by using the DDE-BIFTOOL continuation software  \cite{Engelborghs_01, Engelborghs_02, Heckman}.
In this work we are interested in the details of predicting the appearance of the Hopf bifurcations using approximate perturbation methods.\\  

We offer two derivations of the associated slow flow, one using the two variable expansion perturbation method, and the other by averaging.\\

\begin{figure}[!htbp]
\centering
\includegraphics[width=.6\textwidth]{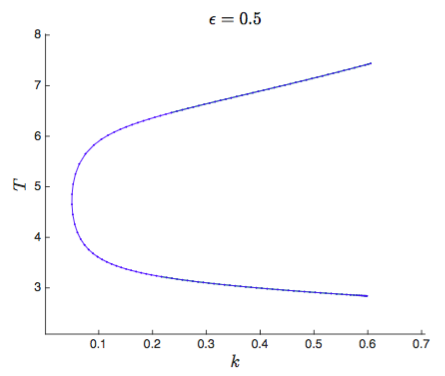}
\caption{Numerical Hopf bifurcation curves for $\epsilon = 0.5$, $\alpha = 0.05$ and $\gamma = 1$ for eq.(\ref{duf}) obtained by using DDE-BIFTOOL .}\label{figure03}
\end{figure}

\section{DERIVATION OF SLOW FLOW}

The two variable method posits that the solution depends on two time variables, $x(\xi,\eta)$, where  $\xi=t$ and $\eta=\epsilon t$.
Then we have 
\begin{equation}
x_d=x(t-T)=x(\xi-T,\eta-\epsilon T)
\end{equation}

Dropping terms of $O(\epsilon^2)$, eq.(\ref{duf}) becomes 
\begin{equation}
x_{\xi\xi} + 2 \epsilon x_{\xi\eta}+x=\epsilon~ \left[-\alpha~x_\xi - \gamma~x^3+k ~x(\xi-T,\eta-\epsilon T)\right]
\end{equation}

Expanding $x$ in a power series in $\epsilon$, $x=x_0+\epsilon x_1 + O(\epsilon^2)$, and collecting terms, we obtain

\begin{equation}
Lx_0 \equiv {x_0}_{\xi\xi}+x_0 = 0
\label{x0eq}
\end{equation}
\begin{equation}
Lx_1 \equiv -2{x_0}_{\xi\eta}-\alpha~{x_0}_\xi - \gamma~{x_0}^3+k ~x_0(\xi-T,\eta-\epsilon T)  
\label{x1eq}
\end{equation}

From eq.(\ref{x0eq}) we have that
\begin{equation}
x_0(\xi,\eta)=A(\eta)\cos\xi+B(\eta)\sin\xi
\label{x0}
\end{equation}
In eq.(\ref{x1eq}) we will need $x_0(\xi-T,\eta-\epsilon T)$:
\begin{equation}
x_0(\xi-T,\eta-\epsilon T)=A_d \cos(\xi-T)+B_d \sin(\xi-T)
\label{x0d}
\end{equation}
where $A_d=A(\eta-\epsilon T)$ and $B_d=B(\eta-\epsilon T)$.

Substituting (\ref{x0}) and (\ref{x0d}) into (\ref{x1eq}) and eliminating resonant terms gives the slow flow:

\begin{eqnarray}
\label{Aeq_duf}
\frac{dA}{d\eta}=-\alpha~\frac{A}{2}+\frac{3~\gamma~B^3}{8}+\frac{\gamma~A^2B}{8}-\frac{k}{2} A_d \sin T-\frac{k}{2} B_d \cos T\\
\label{Beq_duf}
\frac{dB}{d\eta}=-\alpha~\frac{B}{2}-\frac{3~\gamma~A^3}{8}-\frac{\gamma~AB^2}{8}-\frac{k}{2} B_d \sin T+\frac{k}{2} A_d \cos T
\end{eqnarray}

Transforming (\ref{Aeq_duf}),(\ref{Beq_duf}) to polars with $A=R\cos\theta$, $B=R\sin\theta$, we obtain the alternate slow flow:
\begin{eqnarray}
\label{Req}
\frac{dR}{d\eta}&=&-\alpha~\frac{R}{2}-\frac{k}{2} R_d \sin(\theta_d-\theta+ T)\\
\label{thetaeq}
\frac{d\theta}{d\eta}&=&-\frac{3~\gamma~R^2}{2}+\frac{k}{2} \frac{R_d}{R} \cos(\theta_d-\theta+ T)
\end{eqnarray}
where $R_d=R(\eta-\epsilon T)$ and $\theta_d=\theta(\eta-\epsilon T)$.\\

Note that the slow flow (\ref{Req}),(\ref{thetaeq}) contains delay terms in $R_d$ and $\theta_d$ in addition to the usual terms $R$ and $\theta$.
Could this phenomenon be due to some peculiarity of the two variable expansion method?  In order to show that this is not the case, we offer the following slow flow derivation by the method of averaging.\\

We seek a solution to eq.(\ref{duf}) in the form:
\begin{equation}
x(t)=R(t)\cos(t-\theta(t)),~~~\dot x(t)=-R(t)\sin(t-\theta(t))
\label{av1}
\end{equation}
As in the method of variation of parameters, this leads to the (exact) equations:
\begin{eqnarray}
\label{av2}
\frac{dR}{dt}&=&-\epsilon \sin(t-\theta) ~f\\
\label{av3}
\frac{d\theta}{dt}&=&-\frac{\epsilon}{R}\cos(t-\theta)~ f
\end{eqnarray} 
where $f=\alpha~R \sin(t-\theta)-\gamma~R^3\cos^3(t-\theta)+k  R_d\cos(t-T-\theta_d)$,\\
and where  $R_d=R(t-T)$ and $\theta_d=\theta(t-T)$.\\

Now we apply the method of averaging which dictates that we replace the right hand sides of eqs.(\ref{av2}),(\ref{av3}) with averages taken over $2\pi$ in $t$, in which process $R$,$\theta$,$R_d$ and $\theta_d$ are held fixed.  This gives
\begin{eqnarray}
\label{av4}
\frac{dR}{dt}&=&\epsilon\left(-\alpha~\frac{R}{2}-\frac{k}{2} R_d \sin(\theta_d-\theta+ T)\right)\\
\label{av5}
\frac{d\theta}{dt}&=&\epsilon\left( -\frac{3~\gamma~R^2}{2}+\frac{k}{2} \frac{R_d}{R} \cos(\theta_d-\theta+ T)\right)
\end{eqnarray}
Note that  eqs.(\ref{av4}),(\ref{av5}) agree with (\ref{Req}),(\ref{thetaeq}) when $t$ is replaced by $\eta=\epsilon t$.\\

\section{ANALYSIS OF THE SLOW FLOW}

A problem with the slow flow (\ref{Req}),(\ref{thetaeq}) is that they are DDEs rather than ODEs.
Since ODEs are easier to deal with than DDEs, many authors (e.g. Wirkus \cite{Wirkus}, Morrison \cite{Morrison}, Atay \cite{Atay}) simply replace the delay terms by terms with the same variables, but non-delayed.  It is argued that such a step is justified if the product $\epsilon T$ is small:

\begin{eqnarray}
\label{AB}
A_d = A(\eta - \epsilon T) \approx A(\eta) + O(\epsilon),~~~~~~~B_d = B(\eta - \epsilon T) \approx B(\eta) + O(\epsilon).
\end{eqnarray}

\noindent In what follows, we shall refer to this as approach \bf I\rm.
For example, if we replace $A_d$ by $A$, and $B_d$ by $B$, eqs.(\ref{Aeq_duf}),(\ref{Beq_duf}) become:

\begin{eqnarray}
\label{Aeq2}
\frac{dA}{d\eta}=-\alpha~\frac{A}{2}+\frac{3~\gamma~B^3}{8}+\frac{\gamma~A^2B}{8}-\frac{k}{2} A \sin T-\frac{k}{2} B \cos T\\
\label{Beq2}
\frac{dB}{d\eta}=-\alpha~\frac{B}{2}-\frac{3~\gamma~A^3}{8}-\frac{\gamma~AB^2}{8}-\frac{k}{2} B \sin T+\frac{k}{2} A \cos T
\end{eqnarray}

These ODEs have an equilibrium point at the origin.  Linearizing about the origin, we obtain:
\begin{equation}
\frac{d}{d\eta}
\left[
\begin{array}{c}
A\\
B
\end{array}\right]=
\left[
\begin{array}{cc}
-\frac{\alpha}{2}-\frac{k}{2}\sin T&-\frac{k}{2}\cos T\\
\frac{k}{2}\cos T&-\frac{\alpha}{2}-\frac{k}{2}\sin T
\end{array}\right]
\left[
\begin{array}{c}
A\\
B
\end{array}\right]
\label{matrix}
\end{equation}
For a Hopf bifurcation, we require imaginary roots of the characteristic equation, or equivalently (Rand \cite{Rand}, Strogatz \cite{Strogatz}) 
we require the trace of the matrix in eq.(\ref{matrix}) to vanish when the determinant$>0$.
This gives 
\begin{equation}
\mbox{Condition for a Hopf Bifurcation:~~~~~} k \sin T = -\alpha
\label{Hopf_duf}
\end{equation}
Since this condition is based on the bold step of replacing the delay quantities in the slow flow by their undelayed counterparts, the question arises as to the correctness of such a procedure and the validity of eq.(\ref{Hopf_duf}).
See Figure 2 where eq.(\ref{Hopf_duf}) is plotted along with the numerically-obtained conditions for a Hopf.\\

\begin{figure}[!htbp]
\centering
\includegraphics[width=.6\textwidth]{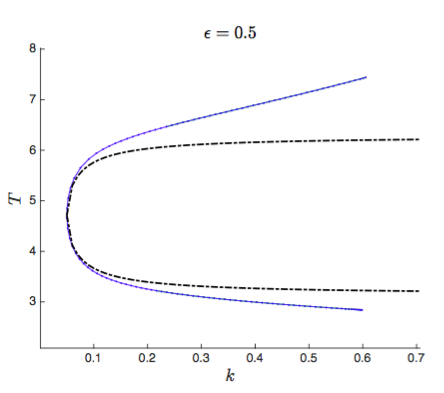}
\caption{Numerical Hopf bifurcation curves (blue/solid) and analytical Hopf condition eq.(\ref{Hopf_duf}) (black/dashdot) for $\epsilon = 0.5$, $\alpha = 0.05$ and $\gamma = 1$ for eq.(\ref{duf}).}\label{figure03}
\end{figure}

Let us now return to eqs.(\ref{Aeq_duf}),(\ref{Beq_duf}) and treat them as DDEs rather than as ODEs.
In what follows we shall refer to this as approach \bf II\rm.
Again linearizing about the origin, we obtain
\begin{eqnarray}
\label{Aeq3}
\frac{dA}{d\eta}=-\frac{\alpha~A}{2}-\frac{k}{2} A_d \sin T-\frac{k}{2} B_d \cos T\\
\label{Beq3}
\frac{dB}{d\eta}=-\frac{\alpha~B}{2}-\frac{k}{2} B_d \sin T+\frac{k}{2} A_d \cos T
\end{eqnarray}
where $A_d=A(\eta-\epsilon T)$ and $B_d=B(\eta-\epsilon T)$.  We set
\begin{equation}
A=a \exp(\lambda \eta),~~~B=b \exp(\lambda \eta),~~~A_d=a \exp(\lambda \eta - \epsilon \lambda T ),~~~B_d=b \exp(\lambda \eta - \epsilon \lambda T )
\label{DDE1}
\end{equation}
where $a$ and $b$ are constants.  This gives
\begin{equation}
\left[
\begin{array}{cc}
-\lambda-\frac{\alpha}{2}-\frac{k}{2}\exp(-\lambda \epsilon T)\sin T&-\frac{k}{2}\exp(-\lambda \epsilon T)\cos T\\
\frac{k}{2}\exp(-\lambda \epsilon T)\cos T&-\lambda-\frac{\alpha}{2}-\frac{k}{2}\exp(-\lambda \epsilon T)\sin T
\end{array}\right]
\left[
\begin{array}{c}
a\\
b
\end{array}\right]=
\left[
\begin{array}{c}
0\\
0
\end{array}\right]
\label{matrix2}
\end{equation}
For a nontrivial solution $(a,b)$ we require the determinant to vanish:
\begin{equation}
\left(-\lambda-\frac{\alpha}{2}-\frac{k}{2}\exp(-\lambda \epsilon T)\sin T \right)^{2}+\frac{k^2}{4}\exp(-2\lambda \epsilon T)\cos^2 T=0
\label{DDE2}
\end{equation}

We set $\lambda=i\omega$ for a Hopf bifurcation and use Euler's formula $\exp(-i\omega \epsilon T)=\cos{\omega\epsilon T}-i \sin{\omega\epsilon T}$. Separating real and imaginary parts we obtain

\begin{eqnarray}
\label{DDE3_duf}
4k^2 \cos 2 \epsilon \omega T+16k\omega \sin T \sin \epsilon \omega T + 8\alpha k \sin T \cos \epsilon \omega T-16\omega^2+4\alpha^2 = 0\\
\label{DDE4_duf}
-4k^2 \sin 2 \epsilon \omega T-8\alpha k \sin T \sin \epsilon \omega T+16k\omega\sin T \cos \epsilon \omega T+16 \alpha \omega=0
\end{eqnarray}

The next task is to analytically solve the two characteristic eqs.(\ref{DDE3_duf})-(\ref{DDE4_duf}) for the pair ($\omega$,$T$). To this aim we use a perturbation schema by setting

\begin{eqnarray}
\label{eq:freq_per_duf}
\omega_{cr} &=&  \sum_{n=0}^N \epsilon^{n} \,\omega_{n}= \omega_{0} + \epsilon\,\omega_{1} + \epsilon^{2}\,\omega_{2} + \dots	\\
\label{eq:delay_per_duf}
T_{cr}&=&  \sum_{n=0}^N \epsilon^{n} \,T_{n} = T_{0} + \epsilon\,T_{1} + \epsilon^{2}\,T_{2} + \dots	
 \end{eqnarray}

Inserting eqs. (\ref{eq:freq_per_duf})-(\ref{eq:delay_per_duf}) in eqs.(\ref{DDE3_duf})-(\ref{DDE4_duf}),\,Taylor expanding the trig functions with respect to the small parameter $\epsilon <<1$, and equating terms of equal order of $\epsilon$ we obtain:

\begin{eqnarray}
\label{eq:freq_exp}
\omega_{cr} &=& \omega_{0} = \frac{\sqrt{k^2-\alpha^2}}{2}	\\
\label{eq:delay_exp}
T_{cr} &=&  T_{0} \left( 1 \pm \epsilon\,\omega_{0} + \epsilon^{2}\,\omega_{0}^{2} \pm \epsilon^{3}\,\omega_{0}^{3} +\epsilon^{4}\,\omega_{0}^{4} \pm\epsilon^{5}\,\omega_{0}^{5}+ \dots	 \right)
 \end{eqnarray}
 
 \noindent where $T_{0}$ is a solution to the equation $\sin T_{0} = -\alpha/k$, that is
 \begin{eqnarray}
\label{eq:delay_T0_duf_01}
T_{0} &=& 2~\pi+\arcsin \left(-\frac{\alpha}{k}\right)\\
\label{eq:delay_T0_duf_02}
T_{0} &=& \pi - \arcsin \left(-\frac{\alpha}{k}\right).
 \end{eqnarray}
(Eqs.(\ref{eq:delay_T0_duf_01})-(\ref{eq:delay_T0_duf_02}) are the black/dashdot curves in Figure 2.)\\

Eq.(\ref{eq:delay_exp}) appears to be the front end of a geometric series.
Assuming the series (\ref{eq:delay_exp}) actually is a geometric series, we can sum it:
 \begin{eqnarray}
\label{eq:delay_series_duf_01}
T_{cr_{1}} =  T_{0} \left( 1 + \epsilon\,\omega_{0} + \epsilon^{2}\,\omega_{0}^{2} + \epsilon^{3}\,\omega_{0}^{3} + \dots	 \right) = \frac{T_{0}}
{1-\epsilon\,\omega_{0}}&&\\
 \nonumber~&&~~|\epsilon\,\omega_{0}| <  1\\
\label{eq:delay_series_duf_02}
T _{cr_{2}}=  T_{0} \left( 1 - \epsilon\,\omega_{0} + \epsilon^{2}\,\omega_{0}^{2} - \epsilon^{3}\,\omega_{0}^{3} + \dots	 \right) = \frac{T_{0}}
{1+\epsilon\,\omega_{0}}
 \end{eqnarray}
 
 Replacing $T_0$ in eqs.(\ref{eq:delay_series_duf_01}),(\ref{eq:delay_series_duf_02}) by the derived values listed in eq.\,(\ref{eq:delay_T0_duf_01})-(\ref{eq:delay_T0_duf_02}), we obtain the following expressions for the critical values $\omega_{cr}$ and $T_{cr}$ for which Hopf bifurcations take place:
 
  \begin{eqnarray}
\label{eq:delay_series_duf_03}
T_{cr_{1}} =   \frac{2\pi+\arcsin \left(-\alpha/k \right)}{1-\epsilon\,\omega_{cr}}&&\\
\nonumber &&~~|\epsilon\,\omega_{cr}| <  1\\
\label{eq:delay_series_duf_04} 
T _{cr_{2}}= \frac{\pi - \arcsin \left(-\alpha/k \right)}{1+\epsilon\,\omega_{cr}}
 \end{eqnarray}
 \noindent where $\omega_{cr} = \omega_{0} = \sqrt{k^2-\alpha^2}/2$.  
Figure 3 shows a comparison of eqs.(\ref{eq:delay_series_duf_03}),(\ref{eq:delay_series_duf_04}) with numerical solutions of eqs.(\ref{DDE3_duf})-(\ref{DDE4_duf})
for various parameters. The numerical solutions were obtained using continuation method. The excellent agreement indicates that eqs.(\ref{eq:delay_series_duf_03}),(\ref{eq:delay_series_duf_04}) are evidently exact solutions of eqs.(\ref{DDE3_duf})-(\ref{DDE4_duf}).\\

\begin{figure}
\centering
\begin{minipage}{8cm}
\includegraphics[width=\textwidth]{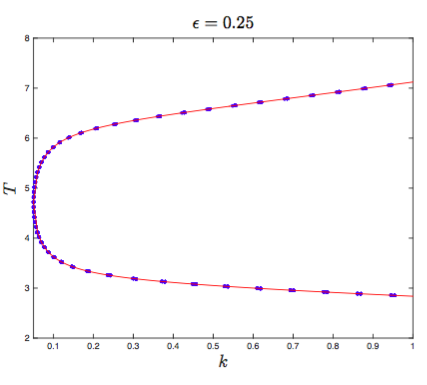}
\end{minipage}
\begin{minipage}{8cm}
\includegraphics[width=\textwidth]{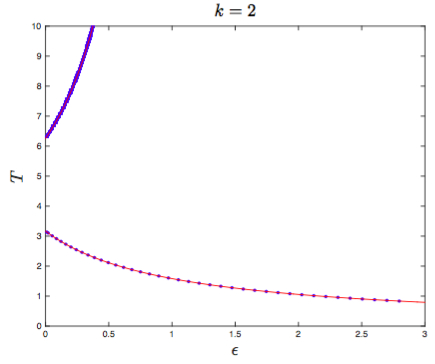}
\end{minipage}
\caption{ Critical delay vs. the feedback magnitude $k$ for $\epsilon = 0.25$ (left). Critical delay vs $\epsilon$ for $k=2$ (right). Red/solid curves: eqs.(\ref{eq:delay_series_duf_03})-(\ref{eq:delay_series_duf_04}). Blue dots: numerical roots of eqs.(\ref{DDE3_duf})-(\ref{DDE4_duf}). These results are for eq.(\ref{duf}) with parameter $\alpha = 0.05$. }\label{figure14}
\end{figure}

\begin{figure}[!htbp]
\centering
\includegraphics[width=.6\textwidth]{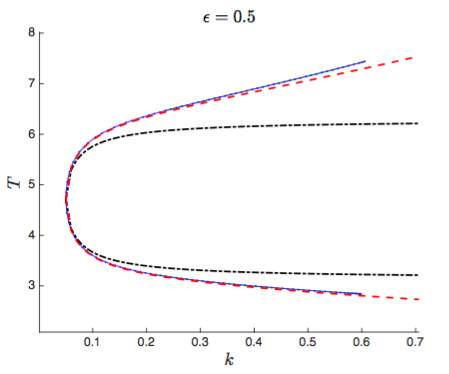}
\caption{Numerical Hopf bifurcation curves (blue/solid) for eq.(\ref{duf}) for $\epsilon = 0.5$, $\alpha = 0.05$ and $\gamma = 1$.  Also shown are the results of approach \bf I\rm, the analytical Hopf condition eq.(\ref{Hopf_duf}) (black/dashdot), and the results of approach \bf II\rm,  eqs.(\ref{eq:delay_series_duf_03}),(\ref{eq:delay_series_duf_04}) (red/dashed) }\label{figure03}
\end{figure}

We now wish to compare the two approaches, namely\\
\bf I \rm: the approach where we replace $A_d$ by $A$, and $B_d$ by $B$ in the slow flow, which gave the condition (\ref{Hopf_duf}), and\\
\bf II \rm: the alternate approach where the terms $A_d$ and $B_d$ are kept without change in the slow flow, 
resulting in eqs.(\ref{eq:delay_series_duf_03}),(\ref{eq:delay_series_duf_04}).\\

Figure 4 shows a comparison between the  analytical Hopf conditions obtained via the two approaches and the numerical Hopf curves. The approach \bf II \rm plotted by red/dashed curves gives a better result than the approach \bf I \rm (black/dashdot curves). Therefore in the case of Duffing equation, treating the slow flow as a DDE gives better results than approximating the DDE slow flow by an ODE. In order to check whether this is also the case for a different type of nonlinearity, we consider in the next section the van der Pol equation with delayed self-feedback.

\section{ANOTHER EXAMPLE: VAN DER POL EQUATION}

As another example we choose the van der Pol equation with delayed self-feedback.
This system has been studied previously by Atay and by Suchorsky et al.
\begin{equation}
\ddot{x}+x=\epsilon ~\left[\dot x(1-x^2)+ k ~x_d\right]
\label{vdp}
\end{equation}

In the case of van der Pol, when there is no feedback ($k=0$), this system is well known to exhibit a stable limit cycle for $\epsilon>0$.  It turns out (Atay \cite{Atay}, Suchorsky \cite{Suchorsky} ) that as delay $T$ increases, for fixed $k>1$, the limit cycle gets smaller and eventually disappears in a Hopf bifurcation. Further increases in $T$ produce another Hopf, which sees the stable limit cycle get reborn.
Figure 5 shows a plot of the Hopfs in the $k-T$ parameter plane. As for the case of Duffing equation we are interested in the details of predicting the appearance of the Hopf bifurcations using approximate perturbation methods. We follow the same procedure as for the case of Duffing equation, that is by deriving the slow flow using the two variable expansion method, and the averaging method. However, for simplicity we omit the averaging method analysis since we obtain the same slow flow by both methods. The obtained slow flow in the cartesian coordinates has the following expression:  

\begin{eqnarray}
\label{Aeq}
\frac{dA}{d\eta}=\frac{A}{2}-\frac{A^3}{8}-\frac{AB^2}{8}-\frac{k}{2} A_d \sin T-\frac{k}{2} B_d \cos T\\
\label{Beq}
\frac{dB}{d\eta}=\frac{B}{2}-\frac{B^3}{8}-\frac{A^2 B}{8}-\frac{k}{2} B_d \sin T+\frac{k}{2} A_d \cos T
\end{eqnarray}

\noindent where $A_d=A(\eta-\epsilon T)$ and $B_d=B(\eta-\epsilon T)$.

\begin{figure}[!htbp]
\centering
\includegraphics[width=.6\textwidth]{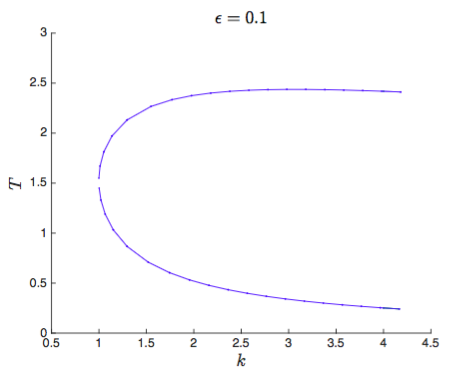}
\caption{Numerical Hopf bifurcation curves for $\epsilon = 0.1$ for eq.(\ref{vdp}) obtained by using DDE-BIFTOOL .}\label{figure03}
\end{figure}

Replacing $A_d$ by $A$, and $B_d$ by $B$, eqs.(\ref{Aeq}),(\ref{Beq}) become:
\begin{eqnarray}
\label{Aeq2}
\frac{dA}{d\eta}=\frac{A}{2}-\frac{A^3}{8}-\frac{AB^2}{8}-\frac{k}{2} A \sin T-\frac{k}{2} B \cos T\\
\label{Beq2}
\frac{dB}{d\eta}=\frac{B}{2}-\frac{B^3}{8}-\frac{A^2 B}{8}-\frac{k}{2} B \sin T+\frac{k}{2} A \cos T
\end{eqnarray}
Linearizing (\ref{Aeq2}) and (\ref{Beq2}) about the origin and looking for the condition where Hopf bifurcation takes place, we find:

\begin{equation}
\mbox{Condition for a Hopf Bifurcation:~~~~~} k \sin T = 1
\label{Hopf}
\end{equation}

This condition is plotted in Figure 6 along with the numerically-obtained conditions for a Hopf.\\

\begin{figure}[!htbp]
\centering
\includegraphics[width=.6\textwidth]{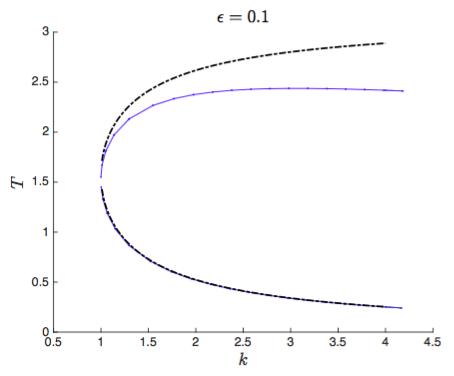}
\caption{Numerical Hopf bifurcation curves (blue/solid) and approach \bf I \rm analytical Hopf condition eq.(\ref{Hopf}) (black/dashdot) for $\epsilon = 0.1$ for eq.(\ref{vdp}) .}\label{figure03}
\end{figure}

If we now treat eqs.(\ref{Aeq}),(\ref{Beq})  as DDEs rather than as ODEs and linearize about the origin, we obtain
\begin{eqnarray}
\label{Aeq3}
\frac{dA}{d\eta}=\frac{A}{2}-\frac{k}{2} A_d \sin T-\frac{k}{2} B_d \cos T\\
\label{Beq3}
\frac{dB}{d\eta}=\frac{B}{2}-\frac{k}{2} B_d \sin T+\frac{k}{2} A_d \cos T
\end{eqnarray}
where $A_d=A(\eta-\epsilon T)$ and $B_d=B(\eta-\epsilon T)$.  We set
\begin{equation}
A=a \exp(\lambda \eta),~~~B=b \exp(\lambda \eta),~~~A_d=a \exp(\lambda \eta - \epsilon \lambda T ),~~~B_d=b \exp(\lambda \eta - \epsilon \lambda T 
\label{DDE1}
\end{equation}
where $a$ and $b$ are constants.  This gives
\begin{equation}
\left[
\begin{array}{cc}
-\lambda+\frac{1}{2}-\frac{k}{2}\exp(-\lambda \epsilon T)\sin T&-\frac{k}{2}\exp(-\lambda \epsilon T)\cos T\\
\frac{k}{2}\exp(-\lambda \epsilon T)\cos T&-\lambda+\frac{1}{2}-\frac{k}{2}\exp(-\lambda \epsilon T)\sin T
\end{array}\right]
\left[
\begin{array}{c}
a\\
b
\end{array}\right]=
\left[
\begin{array}{c}
0\\
0
\end{array}\right]
\label{matrix2}
\end{equation}
For a nontrivial solution $(a,b)$ we require the determinant to vanish:
\begin{equation}
\left(-\lambda+\frac{1}{2}-\frac{k}{2}\exp(-\lambda \epsilon T)\sin T \right)^{2}+\frac{k^2}{4}\exp(-2\lambda \epsilon T)\cos^2 T=0
\label{DDE2}
\end{equation}

We set $\lambda=i\omega$ for a Hopf bifurcation and use Euler's formula $\exp(-i\omega \epsilon T)=\cos{\omega\epsilon T}-i \sin{\omega\epsilon T}$. Separating real and imaginary parts we obtain

\begin{eqnarray}
\label{DDE3}
-\frac{k}{2}\cos\omega\epsilon T\sin T-k\omega \sin\omega\epsilon T\sin T +\frac{k^2}{4}\cos{2\omega\epsilon T}+\frac{1}{4}-\omega^2=0\\
\label{DDE4}
k\omega \cos\omega\epsilon T\sin T 
+\frac{k}{2}\sin\omega\epsilon T\sin T-\frac{k^2}{4}\sin{2\omega\epsilon T}-\omega=0
\end{eqnarray}

As in the case of Duffing equation we proceed by using a perturbation schema to analytically solve the two characteristic eqs.\,(\ref{DDE3})-(\ref{DDE4}) for the pair ($\omega$,$T$). We set the critical frequency and delay to be: 

\begin{eqnarray}
\label{eq:freq_per}
\omega_{cr} &=&  \sum_{n=0}^N \epsilon^{n} \,\omega_{n}= \omega_{0} + \epsilon\,\omega_{1} + \epsilon^{2}\,\omega_{2} + \dots	\\
\label{eq:delay_per}
T_{cr}&=&  \sum_{n=0}^N \epsilon^{n} \,T_{n} = T_{0} + \epsilon\,T_{1} + \epsilon^{2}\,T_{2} + \dots	
 \end{eqnarray}
 
  \noindent where $T_{0}$ is a solution to the equation $\sin T_{0} = 1/k$, that is
 \begin{eqnarray}
\label{eq:delay_T0_01}
T_{0} &=& \arcsin \left(\frac{1}{k}\right)\\
\label{eq:delay_T0_02}
T_{0} &=& \pi - \arcsin \left(\frac{1}{k}\right).
 \end{eqnarray}
(Eqs.(\ref{eq:delay_T0_01})-(\ref{eq:delay_T0_02}) are the black/dashdot curves in Figure 6.)\\

Inserting eqs. (\ref{eq:freq_per})-(\ref{eq:delay_per}) in eqs.(\ref{DDE3})-(\ref{DDE4}),\,Taylor expanding the trig functions with respect to the small parameter $\epsilon <<1$, and equating terms of equal order of $\epsilon$ we obtain:

\begin{eqnarray}
 \label{eq:freq_exp}
\omega_{cr} &=& \omega_{0} = \frac{\sqrt{k^2-1}}{2}	\\
\label{eq:delay_series_03}
T_{cr_{1}} &=&   \frac{\arcsin \left(1/k \right)}{1-\epsilon\,\omega_{cr}}\\
&&~~~~~~~~~~~~~~~~~~~~~~~~~~~~~~|\epsilon\,\omega_{cr}| <  1\nonumber \\
\label{eq:delay_series_04} 
T _{cr_{2}}&=& \frac{\pi - \arcsin \left(1/k \right)}{1+\epsilon\,\omega_{cr}}
 \end{eqnarray}

Figure 7 shows a comparison of eqs.(\ref{eq:delay_series_03}),(\ref{eq:delay_series_04}) with numerical solutions of eqs.(\ref{DDE3})-(\ref{DDE4}) for various parameters. The excellent agreement indicates that eqs.(\ref{eq:delay_series_03}),(\ref{eq:delay_series_04}) are evidently exact solutions of eqs.(\ref{DDE3})-(\ref{DDE4}).\\

\begin{figure}
\centering
\begin{minipage}{8cm}
\includegraphics[width=\textwidth]{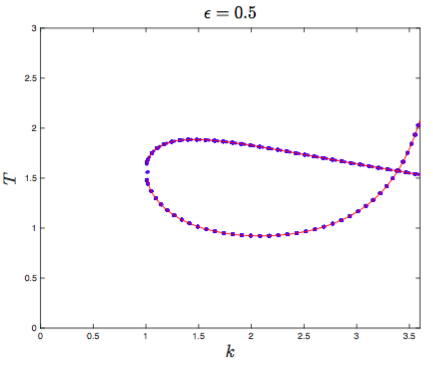}
\end{minipage}
\begin{minipage}{8cm}
\includegraphics[width=\textwidth]{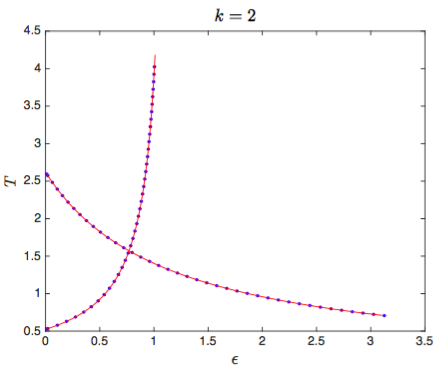}
\end{minipage}
\caption{ Comparison of numerical versus analytic results obtained by approach \bf II \rm for eq.(\ref{vdp}). Critical delay vs. the feedback magnitude $k$ for $\epsilon = 0.5$ (left). Critical delay vs $\epsilon$ for $k=2$ (right). Red/solid curves: eq.(\ref{eq:delay_series_03})-(\ref{eq:delay_series_04}). Blue dots: numerical roots of eqs.(\ref{DDE3})-(\ref{DDE4}). }\label{figure14}
\end{figure}

\begin{figure}
\begin{minipage}{8cm}
\includegraphics[width=\textwidth]{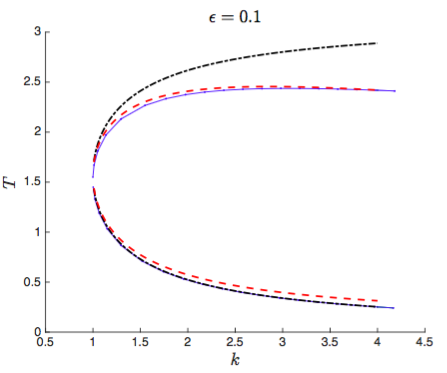}
\end{minipage}
\begin{minipage}{8cm}
\includegraphics[width=\textwidth]{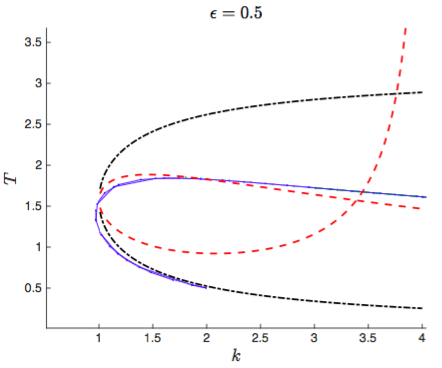}
\end{minipage}
\caption{Numerical Hopf bifurcation curves (blue/solid) for eq.(\ref{vdp}), for $\epsilon$= 0.1 (left) and $\epsilon$= 0.5 (right).  Also shown are the results of approach \bf I\rm, the analytical Hopf condition eq.(\ref{Hopf}) (black/dashdot), and the results of approach \bf II\rm,  eqs.(\ref{eq:delay_series_03})-(\ref{eq:delay_series_04}) (red/dashed).}\label{figure14x}
\end{figure}

Figure 8 shows a comparison between the numerically-obtained Hopf conditions and the Hopf conditions obtained by following the two approaches \bf I \rm and \bf II\rm. When $\epsilon = 0.1$, eq.\,(\ref{eq:delay_T0_01}) (black/dashdot curve) gives a perfect match with the lower numerical branch (blue/solid curve)  than eq.\,(\ref{eq:delay_series_03}) (red/dashed curve). However for the upper numerical branch, eq.\,(\ref{eq:delay_series_04}) gives a better approximation than eq.\,(\ref{eq:delay_T0_02}), see Figure 8. As $\epsilon$ is increased ($\epsilon = 0.5$),  eq.\,(\ref{eq:delay_T0_01}) still gives a better approximation for the lower numerical branch than eq.\,(\ref{eq:delay_series_03}). On the other hand eq.\,(\ref{eq:delay_series_04}) succeeds in tracking the upper numerical branch, see Figure 8.

\section{DISCUSSION }

In the two studied examples we saw that the two approaches gave different results. In the Duffing equation, the approach \bf II \rm gave better results.  This is expected since we did not approximate $A_d$ by $A$, and $B_d$ by $B$, and instead analyzed the slow flow as a DDE. However in the van der Pol example, we obtained unexpected results. From Figure 8, the upper Hopf branch obtained by the approach \bf II \rm gave a better approximation of the upper numerical Hopf curve than the one obtained from approach \bf I\rm, eq.\,(\ref{eq:delay_T0_02}). This could be explained by the fact that as $T$ is increased, the term $\epsilon\,T$ increases as well, which makes the approximation  $A_d = A$, and $B_d = B$ no longer valid. By contrast, in the the approach \bf II \rm  the increasing of $T$ does not affect  the condition (\ref{eq:delay_series_04}). But unexpectedly, the condition (\ref{eq:delay_series_03}) obtained by approach \bf II \rm fails to give a better result for the lower Hopf curve. This could be explained by the singularity that takes place in the lower numerical Hopf branch where the limit cycle disappears and the origin $x=0$ changes its nature as an equilibrium. For example when $\epsilon = 0.5$ this singular behavior occurs for $k \geq 2$, see Figure 8.  Both the method of averaging and the two variable expansion perturbation method are built on the assumption that the solution at $O(\epsilon^{0})$ is a periodic solution around the origin $x=0$.  However for increasing $\epsilon$ and $k$ the origin no longer exhibits this behavior, and our assumption of the periodicity of our unperturbed solution does not hold anymore. Note that eq.\,(\ref{eq:delay_T0_01}) does not contain an $\epsilon$ term, thus it does not vary with increasing $\epsilon$. \\

\noindent Figure 9 shows a numerical simulation of the van der Pol equation (\ref{vdp}) for $k=2.1$, where the origin has changed its nature. This figure corresponds to the lower Hopf curve in Figure 8 when $\epsilon = 0.5$. This unexpected failure of approach \bf II \rm leads us to wonder if this happens because the system is a self-sustained one. In order to show that is not the case, we consider a limit cycle system studied by Erneux and Grasman \cite{Erneux}. In their work, they looked for the Hopf curves in a limit cycle system with delayed self-feedback:

\begin{equation}
\ddot{x}+x=\epsilon ~\left[\dot x(1-x^2)+ k ~x_d - k ~x \right]
\label{erneux}
\end{equation}

We apply the same procedure, approach \bf II\rm, as we did for Duffing and van der Pol examples to equation eq.\,(\ref{erneux}), and we obtain the following critical frequency and time delay:

 \begin{eqnarray}
\label{eq:delay_w1_erneux}
\omega_{cr_{1}}&=& \sqrt{ \frac{k^2}{2} -\frac{1}{4} - \frac{k}{2} \sqrt{k^2-1} } \\
\label{eq:delay_w2_erneux}
\omega_{cr_{2}}&=& \sqrt{ \frac{k^2}{2} -\frac{1}{4} + \frac{k}{2} \sqrt{k^2-1} }
 \end{eqnarray}
 
\begin{eqnarray}
\label{eq:delay_T1_erneux}
T _{cr_{1}}&=& \frac{\pi - \arcsin \left(1/k \right)}{1+\epsilon\,\omega_{cr_{1}}}\\
\label{eq:delay_T2_erneux}
T _{cr_{1}}&=& \frac{ \arcsin \left(1/k \right)}{1+\epsilon\,\omega_{cr_{2}}}
 \end{eqnarray}

\noindent  Figure 10 shows a comparison between approach \bf II\rm, eqs.\,(\ref{eq:delay_w1_erneux}),(\ref{eq:delay_w2_erneux}),(\ref{eq:delay_T1_erneux}),(\ref{eq:delay_T2_erneux}), and approach \bf I\rm, which again gives eqs.\,(\ref{eq:delay_T0_01}),(\ref{eq:delay_T0_02}), and the numerical Hopf curves obtained by use of DDE-BIFTOOL. Figure 10 shows that approach \bf II \rm gives better results than approach \bf I\rm. However approach \bf I \rm still gives a good fit for the lower Hopf curve as in the case of eq.\,(\ref{vdp}).

\begin{figure}[!htbp]
\centering
\includegraphics[width=.6\textwidth]{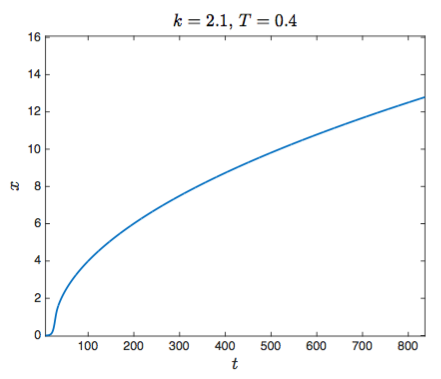}
\caption{Numerical integration for $x$ as a function of time $t$ in eq.(\ref{vdp}) for $\epsilon = 0.5$,  $k=0.21$ and delay $T$=0.4.  Note that the motion grows large and there is no limit cycle.  The origin has changed its nature. See Figure 8 and text. }\label{figure03}
\end{figure}

\begin{figure}
\begin{minipage}{8cm}
\includegraphics[width=\textwidth]{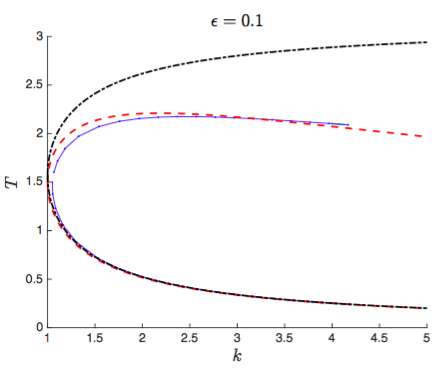}
\end{minipage}
\begin{minipage}{8cm}
\includegraphics[width=\textwidth]{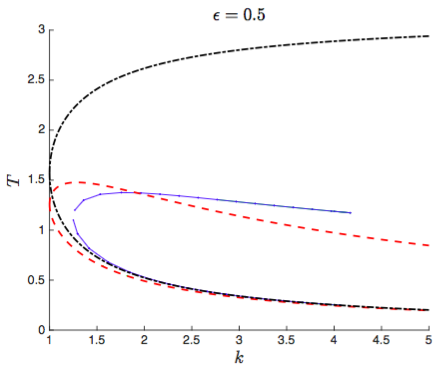}
\end{minipage}
\caption{Numerical Hopf bifurcation curves (blue/solid) for eq.(\ref{erneux}) for $\epsilon$= 0.1 (left) and $\epsilon$= 0.5 (right).  Also shown are the results of approach \bf I\rm, the analytical Hopf condition eq.(\ref{Hopf}) (black/dashdot), and the results of approach \bf II\rm,   eqs.(\ref{eq:delay_T1_erneux})-(\ref{eq:delay_T2_erneux}) (red/dashed).}\label{figure14}
\end{figure}

\section{CONCLUSION	}

When a DDE with delayed self-feedack is treated using a perturbation method (such as the two variable expansion method, multiple scales, or averaging), the resulting slow flow typically involves delayed variables.  In this work we compared the behavior of the resulting DDE slow flow with a related ODE slow flow obtained by  replacing the delayed variables in the slow flow with non-delayed variables.  We studied sample systems based on the Duffing equation with delayed self-feedback, eq.(\ref{duf}), and on the van der Pol equation with delayed self-feedback, eq.(\ref{erneux}).  In both cases we found that replacing the delayed variables in the slow flow by non-delayed variables (approach \bf I \rm) gave better results on the lower Hopf curve than on the upper Hopf curve.\\

Our conclusion is therefore that the researcher is advised to perform the more lengthy approach \bf II \rm analysis on the DDE slow flow in situations where values of the product $\epsilon T$ is relatively large, as in the upper Hopf curves in Figures 1.


\begin{thebibliography}{00}



\bibitem{Ji}  J.C. Ji, A.Y.T. Leung, Resonances of a nonlinear SDOF system with two time-delays on linear feedback control, Journal of Sound and
Vibration 253,  985-1000 (2002)

\bibitem{Maccari} A. Maccari, The resonances of a parametrically excited Van der Pol oscillator to a time delay state feedback,
Nonlinear Dynamics 26, 105-119 (2001)

\bibitem{Hu} Hu, H., Dowell, E. H., and Virgin, L. N., Resonances of a harmonically forced duffing oscillator with time delay state
feedback, Nonlinear Dynamics 15, 311-327 (1998) 

\bibitem{Wahi}  Wahi, P., Chatterjee, A., Averaging oscillations with small fractional damping and delayed terms. Nonlinear Dyn. 38, 3-22 (2004)

\bibitem{Atay} Atay, F.M., Van der Pol's oscillator under delayed feedback. J. Sound Vib. 218(2), 333-339 (1998)

\bibitem{Suchorsky} Suchorsky, M.K., Sah, S.M., Rand, R.H.,Using delay to quench undesirable vibrations. Nonlinear Dyn. 62, 107-116 (2010)


 \bibitem{Engelborghs_01}  Engelborghs, K., Luzyanina, T., Roose, D., Numerical bifurcation analysis of delay differential equations using: DDE- BIFTOOL. ACM Trans. Math. Softw. 28(1), 1-21 (2002)
 
 
 
 \bibitem{Engelborghs_02}  Engelborghs, K., Luzyanina, T., Samaey, G., DDE-BIFTOOL v. 2.00: a Matlab package for bifurcation analysis of delay differential equations. Technical Report TW-330, Dept. Comp. Sci., K.U.Leuven, Leuven, Belgium (2001)
 
 
 
 \bibitem{Heckman}   Heckman, C.R., An introduction to DDE-BIFTOOL is available as Appendix B of the doctoral thesis of Christoffer Heckman: asymptotic and numerical analysis of delay- coupled microbubble oscillators (Doctoral Thesis). Cornell University (2012)

\bibitem{Wirkus} Wirkus, S., Rand, R.H., The dynamics of two coupled van der Pol oscillators with delay coupling. Nonlinear Dyn. 30, 205-221 (2002)


\bibitem{Morrison} Morrison, T.M., Rand, R.H., 2:1 Resonance in the delayed nonlinear Mathieu equation. Nonlinear Dyn. 50, 341-352 (2007)

\bibitem{Rand} Rand, R.H.,  Lecture notes in nonlinear vibrations pub- lished on-line by the Internet-First University Press http:// ecommons.library.cornell.edu/handle/1813/28989 (2012)

\bibitem{Strogatz} Strogatz, S. H., Nonlinear Dynamics and Chaos: With Applications to Physics, Biology, Chemistry, and Engineering (Addison-Wesley, Reading, Mas- sachusetts) (1994)


\bibitem{Moon} Kalmar-Nagy, T., Stepan, G. and Moon, F.C.,  Subcritical Hopf bifurcation in the delay equation model for machine tool vibrations,
Nonlinear Dynamics 26:121-142 (2001)

\bibitem{Verdugo} Verdugo, A. and Rand, R., Hopf Bifurcation in a DDE Model of Gene Expression, Communications in Nonlinear Science and Numerical Simulation 13:235-242 (2008)


\bibitem{Erneux}  T. Erneux, J. Grasman, Limit-cycle oscillators subject to a delayed feedback, Phys. Rev. E 78 (2)  026209 (2008)

 
\end{thebibliography}
\end{document}